\documentclass[12pt]{article}
\usepackage{amscd}
\usepackage{amsmath,amsfonts,amssymb,amscd}
\usepackage{indentfirst,graphicx,epsfig}
\usepackage{graphicx,psfrag}
\makeatletter \@addtoreset{figure}{} \makeatother \makeatletter
\long\def\@makecaption#1#2{%
   \vskip 10\p@
   \setbox\@tempboxa\hbox{{#1}\ \ #2}%
   \ifdim \wd\@tempboxa >\hsize
       {#1}\ \ #2\par
   \else
       \hbox to\hsize{\hfil\box\@tempboxa\hfil}%
   \fi}
\makeatother

\newtheorem{thm}{Theorem}
\newtheorem{pro}{Proposition}
\newtheorem{cor}{Corollary}
\newtheorem{lem}{Lemma}

\makeatletter \@addtoreset{equation}{section}

\def\qed{\hfill \rule{4pt}{7pt}}
\def\pf{\noindent {\it Proof.} }
\begin{document}

\begin{center} {\Large \bf Note on the $4$- and $5$-leaf
powers\footnote{Supported by NSFC No.10831001, PCSIRT and the
``973" program. }}
\end{center}
\pagestyle{plain}
\begin{center}
{
  {\small Xueliang Li, Yongtang Shi, Wenli Zhou}\\

  {\small Center for Combinatorics and LPMC-TJKLC}\\
  {\small Nankai University, Tianjin 300071, China}\\
  {\small  Email: lxl@nankai.edu.cn}\\
    }
\end{center}

\begin{center}
\begin{minipage}{120mm}
\vskip 0.3cm
\begin{center}
{\bf Abstract}
\end{center}
{\small Motivated by the problem of reconstructing evolutionary
history, Nishimura et al. defined $k$-leaf powers as the class of
graphs $G=(V,E)$ which has a $k$-leaf root $T$, i.e., $T$ is a tree
such that the vertices of $G$ are exactly the leaves of $T$ and two
vertices in $V$ are adjacent in $G$ if and only if their distance in
$T$ is at most $k$. It is known that leaf powers are chordal graphs.
Brandst\"adt and Le proved that every $k$-leaf power is a
$(k+2)$-leaf power and every $3$-leaf power is a $k$-leaf power for
$k\geq 3$. They asked whether a $k$-leaf power is also a
$(k+1)$-leaf power for any $k\geq 4$. Fellows et al. gave an example
of a $4$-leaf power which is not a $5$-leaf power. It is interesting
to find all the graphs which have both $4$-leaf roots and $5$-leaf
roots. In this paper, we prove that, if $G$ is a $4$-leaf power with
$L(G)\neq \emptyset$, then $G$ is also a $5$-leaf power, where
$L(G)$ denotes the set of leaves of $G$.
 }\\[3mm]

\end{minipage}
\end{center}

\section{Introduction}

Motivated by the problem of reconstructing evolutionary history,
Nishimura et al. \cite{NRT} introduced the notion of $k$-leaf root
and $k$-leaf power. Let $G=(V(G), E(G))$ be a finite undirected
graph. For $k\geq 2$, a tree $T$ is a {\it $k$-leaf root} of $G$
if $V(G)$ is the set of leaves of $T$ and two vertices $x, y\in
V(G)$ are adjacent in $G$ if and only if their distance $d_T(x,
y)$ in $T$ is at most $k$, i.e., $xy\in E(G)$ if and only if
$d_T(x, y)\leq k$. $G$ is called the {\it $k$-leaf power} of $T$.

Obviously, a graph $G$ is a $2$-leaf power if and only if it is the
disjoint union of cliques, i.e., $G$ is $P_3$-free. The $3$-leaf
powers are exactly the bull-, dart-, and gem-free chordal graphs
\cite{DGHN}; equivalently, $3$-leaf powers are exactly the result of
substituting cliques into the nodes of a tree. For more results on
$3$-leaf powers see \cite{BL, R}. A characterization of 4-leaf
powers in terms of forbidden subgraphs is much more complicated
\cite{BLS0, R}. For $5$-leaf powers, a polynomial time recognition
was given in \cite{CK} but no structural characterization of
$5$-leaf powers has been known. Recently, the authors of \cite{BLR}
characterized distance-hereditary $5$-leaf powers (without pair of
similar vertices) in terms of forbidden induced subgraphs. The
complexity of characterizing and recognizing leaf powers in general
is a major open problem, and so the complexity of characterizing and
recognizing leaf powers for $k\geq 6$ is a major open problem.

It is well known that if a graph $G$ has $k$-leaf roots for some
$k$, then $G$ is chordal. However, there are chordal graphs which
are not $k$-leaf power for any $k\geq 2$. In \cite{BH}, the authors
showed that every interval graph is a leaf power and unit interval
graphs are exactly the leaf powers which have a caterpillar as leaf
root.

Brandst\"adt and Le \cite{BL} proved that every $k$-leaf power is a
$(k+2)$-leaf power and every $3$-leaf power is is a $k$-leaf power
for $k\geq 3$. They asked whether $k$-leaf power is also a
$(k+1)$-leaf power for any $k\geq 4$. Fellows et al. gave an example
of a $4$-leaf power which is not a $5$-leaf power. Then it is
interesting to find all the graphs which have both $4$-leaf roots
and $5$-leaf roots. In this paper, we prove that, if $G$ is a
$4$-leaf power with $L(G)\neq \emptyset$, then $G$ is also a
$5$-leaf power, where $L(G)$ denotes the set of leaves of $G$.

\section{Notations and basic facts}

We consider $G=(V(G), E(G))$ as a finite, simple and undirected
graph. For $k\geq 1$, let $P_k$ denote a path with $k$ vertices and
$k-1$ edges, and, for $k\geq 3$, let $C_k$ denote a cycle with $k$
vertices and $k$ edges. A vertex $v$ is {\it pendent} or a {\it
leaf} if its degree is $1$, i.e., $d(v)=1$. Let $L(G)$ denotes the
set of leaves of $G$.  A path $v_0v_1\ldots v_t$ ($t\geq 1$) is {\it
pendent} if $d(v_0)=1$, $d(v_t)\geq 3$ and all other vertices have
degree $2$. Especially, if $t=1$, it is called a {\it pendent edge}.
The neighborhood of a vertex $u\in V(G)$ in the graph $G$ is denoted
by $N_G(u)$ and the closed neighborhood is denoted by
$N_G[u]=\{u\}\cup N_G(u)$. Two vertices $u,v\in V(G)$ with $u\neq v$
are called {\it similar} if $N_G[u]=N_G[v]$. It is obvious that if
$u$ and $v$ are a pair of similar vertices, then $uv\in E(G)$. An
undirected graph is {\it chordal (triangulated, rigid circuit)} if
every cycle of length greater than three has a chord, which is an
edge connecting two nonconsecutive vertices on the cycle. Namely, a
graph is chordal if it contains no induced $C_k$ for $k\geq 4$. A
vertex is {\it simplicial} in $G$ if its closed neighborhood $N[v]$
is a clique. It is well known that every chordal graph $G$ has a
simplicial vertex, furthermore, if $G$ is not complete, then it has
two nonadjacent simplicial vertices \cite{D}. For various
characterizations of chordal graphs, we refer to \cite{BLS}.

Firstly, we list some useful facts on leaf powers in Proposition
\ref{pro1} without proofs.

\begin{pro}(\cite{BL})\label{pro1} (i) \ For every $k\geq 2$, $k$-leaf powers are
chordal.\\
(ii) \ Every induced subgraph of a $k$-leaf power is a $k$-leaf
power for $k\geq 2$.\\
(iii) \ A graph is a $k$-leaf power if and only if each of its
connected components is a $k$-leaf power.
\end{pro}

We mention some basic facts from \cite{BL, R}, and repeat their
proofs for completeness.

\begin{pro}(\cite{BL})\label{pro2} (i) \ Every $k$-leaf power is a $(k+2)$-leaf power.\\
(ii) \ Every $3$-leaf power is a $k$-leaf power for $k\geq 3$.
\end{pro}
\pf (i) \ Let $T$ be a $k$-leaf root of $G$, and let $T'$ be the
tree obtained from $T$ by subdividing each pendant edge with a new
vertex. Thus, the leaves of $T'$ are exactly those of $T$. Clearly,
for all $x,y\in V(G)$, $xy\in E(G)$ if and only if $k\geq d_T(x, y)
=d_{T'}(x,y)-2$, hence $T'$ is a $(k+2)$-leaf of $G$.

(ii) \ Let $T$ be a $3$-leaf root of a graph $G$, and let $T'$ be
the tree obtained from $T$ by subdividing each non-pendent edge with
exactly $k-3$ new vertices. Thus, the leaves of $T'$ are exactly
those of $T$. Clearly, for all $x,y\in V(G)$, $xy\in E(G)$ if and
only if $d_T(x,y)= d_{T'}(x,y)=2$, or $d_T(x,y)=3$ and $d_{T'}(x,y)
=k$. Hence $T'$ is a $k$-leaf root of $G$.\qed

\begin{pro}(\cite{R})\label{pro3}
If $G=(V,E)$ is a graph and $u,v\in V$ are similar, then $G$ has a
$k$-leaf root if and only if $G-u$ has a $k$-leaf root for any
$k\geq 2$.
\end{pro}
\pf If $G$ has a $k$-leaf root $T$, denoting $T'$ the tree obtained
from $T$ by deleting the leaf $u$, i.e., the pendent path containing
$u$ in $T$, then $T'$ is a $k$-leaf root of $G-u$.

If $G-u$ has a $k$-leaf root $T'$, we can obtained a leat root $T$
of $G$ by attaching the new leaf $u$ to the neighbor of the leaf $v$
in $T'$.\qed

\begin{figure}[ht]
\begin{center}
\includegraphics[width=14cm]{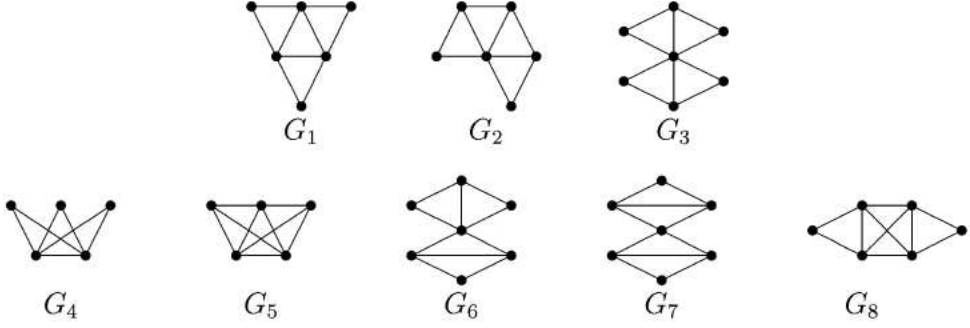}
\end{center}
\caption{Forbidden subgraphs of $4$-leaf powers without similar
vertices.} \label{fig1}
\end{figure}

\begin{thm}[\cite{R}]\label{thm0}
Let $G=(V,E)$ be a graph without pairs of similar vertices. Then $G$
has a $4$-leaf root if and only if it is chordal and does not
contain any of the graphs in Figure \ref{fig1} as an induced
subgraph.
\end{thm}

\section{Main results}

\begin{figure}[ht]
\psfrag{u}{$u$}\psfrag{c1}{$c_1$}\psfrag{c2}{$c_2$}
\psfrag{ct}{$c_t$}\psfrag{v}{$v$}
\psfrag{H}{$H$}\psfrag{T0}{$T_0$}\psfrag{C}{$C$}
\begin{center}
\includegraphics[width=10cm]{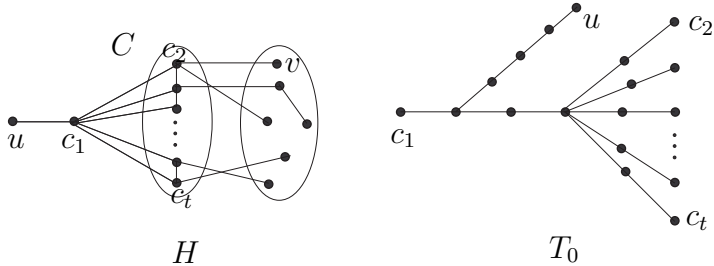}
\end{center}

\caption{One $5$-leaf root of a graph $H\in \mathcal{H}$.}
\label{fig2}
\end{figure}

We define a class of simple graphs $\mathcal{H}$ satisfying that:
$H\in \mathcal{H}$, if $H$ has a clique $C$ with vertex set
$\{c_1,c_2,\ldots,c_t\}$ such that there is a leaf $u$ adjacent to
one vertex $c_1$ of $C$ and $|N_H(v)\cap C|\leq 1$ for all $v\in
H\setminus \{u,c_1,\ldots,c_t\}$, $t\geq 2$. Let $T_0$ be a tree (as
shown in Figure \ref{fig2}) satisfying that $d_{T_0}(c_1,u)=5$,
$d_{T_0}(c_1,c_i)=5$ for $2\leq i\leq t$ and $d_{T_0}(c_i,c_j)=4$
for all $2\leq i,j\leq t$ and $i\neq j$ if $t\geq 3$.

\begin{lem}\label{lem1}
For any graph $H\in \mathcal{H}$, suppose $H$ is a $5$-leaf power
without pairs of similar vertices, then there must be a $5$-leaf
root $T$ of $H$ containing $T_0$ as its subtree.
\end{lem}
\pf Let $H\in \mathcal{H}$ be a $5$-leaf power without pairs of
similar vertices. We can assume that in each $5$-leaf root $T$ of
$H$, the pendent path ending at the leaf $u$ must have length $4$
and connect to the neighbor of $c_1$ in $T$. Otherwise, since $u$ is
a leaf of $H$, we can attain what we need by subdividing the leaf
edge and moving the pendent paths ending at the leaf $u$ to the
neighbor of $c_1$. Thus, suppose $T$ is a $5$-leaf root of $H$
satisfying the above conditions. If $T$ contains $T_0$ as a subtree,
we are done. Otherwise, we will construct a new tree $T'$ from $T$,
which is a $5$-leaf root of $H$ and contains $T_0$ as a subtree.

If $t=2$, we have $d_T(c_1,c_2)\geq 3$. Since $T$ does not contain
$T_0$ as a subtree, then $d_T(c_1,c_2)=3$ or $4$. Since $H$ have no
pairs of similar vertices, then in these two cases, the pendent path
containing $u$ in $T$ must connect to the neighbor of $c_1$. Let
$c_1abc_2$ (or $c_1abcc_2$) be the path joining vertices $c_1$ and
$c_2$ in $T$. Denote by $T'$ the new tree constructed from $T$ by
subdividing edge $ab$ with exactly two (or one) new vertices. It is
easy to verify that $T'$ is a $5$-leaf root of $H$, which contains
$T_0$ as a subtree. Thus, in the following we always assume $t\geq
3$.

If for each $2\leq i\leq t$, $d_T(c_1,c_i)=3$, without loss of
generality, let $c_1abc_2$ be the path joining vertices $c_1$ and
$c_2$ in $T$. Since $H$ has no pairs of similar vertices, we notice
that in $T$ there are no pendent paths containing $c_i$ ($i\geq 3$)
and connecting to the vertex $b$ and the length of each pendent path
containing $c_i$ ($i\geq 3$) and connecting to the vertex $a$ must
be $2$. Denote by $T'$ the new tree constructed from $T$ by
subdividing edge $c_1a$ with exactly two new vertices, deleting the
pendent path containing $u$ and attaching a pendent path with leaf
$u$ to the neighbor of $c_1$ such that $d_{T'}(c_1,u)=5$. It is easy
to verify that $d_{T'}(c_1,u)=5$, $d_{T'}(c_1,c_i)=5$ for $2\leq
i\leq t$ and $d_{T'}(c_i,c_j)=4$ for all $2\leq i,j\leq t$ and
$i\neq j$. Therefore, $T'$ is a $5$-leaf root of $H$, which contains
$T_0$ as a subtree.

Assume $\max_{2\leq i\leq t} d_T(c_1, c_i)=4$, without loss of
generality, suppose $d_T(c_1, c_2)=4$. Let $c_1abcc_2$ be the path
joining vertices $c_1$ and $c_2$ in $T$, then there are no pendent
paths containing $c_i$ ($i\geq 3$) and connecting to the vertex $c$.
Let $P_{ai}$ (or $P_{bi}$) be the pendent paths containing $c_i$ and
connecting to the vertex $a$ (or $b$). Since $\max_{2\leq i\leq t}
d_T(c_1, c_i)=4$ and $C$ is a clique, the length of $P_{ai}$ (or
$P_{bi}$) is at most $2$. Since $H$ has no pairs of similar
vertices, the length of each $P_{ai}$ is exactly two and there is at
most one pendent path with length one among all pendent paths
$P_{bi}$ ($i\geq 3$). Without loss of generality, let $|P_{b3}|=1$.
Denote by $T'$ the new tree obtained from $T$ by subdividing edge
$bc_3$ with exactly one new vertex, then moving all the pendent
paths $P_{ai}$ ($i\geq 4$, if they exist) to the vertex $b$, and
then subdividing edge $ab$ with exactly one new vertex. We can
verify that $T'$ is a $5$-leaf root of $H$, which contains $T_0$ as
a subtree. Firstly, we notice that for all $2\leq i,j\leq t$ and
$i\neq j$, $d_{T'}(c_1,c_i)=5$, $d_{T'}(c_i,c_j)=4$. For each vertex
$v\in H\setminus \{u,c_1,\ldots,c_t\}$ such that $d_{T'}(v,c_i)\leq
5$, there exists no vertex $c_j$ ($j\neq i$) such that
$d_{T'}(v,c_j)\leq 5$. Otherwise, we have $d_T(v,c_1)\leq 5$ or
$d_{T}(v,c_j)\leq 5$, contradicting to the definition of
$\mathcal{H}$.

For the case $\max_{2\leq i\leq t} d_T(c_1, c_i)=5$, we also can
construct a new tree $T'$ from $T$ by the similar transformations as
above, which contains $T_0$ as a subtree. \qed

\begin{thm}\label{thm1}
If $G$ is a $4$-leaf power with $L(G)\neq \emptyset$, then $G$ is
also a $5$-leaf power, where $L(G)$ denotes the set of leaves of
$G$.
\end{thm}
\pf In view of Proposition \ref{pro1} (iii) and Proposition
\ref{pro3}, it suffices to consider connected graphs without pairs
of similar vertices.

By contrary, we assume $G=(V,E)$ is a $4$-leaf power without pairs
of similar vertices such that $L(G)\neq \emptyset$, and $G$ is not a
$5$-leaf power. Furthermore, we assume that among all such graphs
$G$ is chosen such that it has minimum number of edges and subject
to this condition it has maximum number of leaves.

Let $u$ be a leaf of $G$ and $v_1$ the unique neighbor of $u$. If
$G'=G-u$ has no pair of similar vertices, then $G'$ has a $5$-leaf
root $T'$. Thus, we can obtain a $5$-leaf root of $G$ by attaching a
pendent path of length $4$ ending at the new leaf $u$ to the
neighbor of the leaf $v_1$ in $T'$.

Hence, we may assume that $G'$ has a pair of similar vertices. Note
that $v_1$ must be one of each pair of similar vertices, since $v_1$
is the only vertex whose neighborhood is changed in $G'$. In fact,
there is a unique vertex $v_2\in V$ such that $v_1$ and $v_2$ are
similar in $G'$. For otherwise, suppose $v_1$ and $v_3$ ($v_3\neq
v_2$) are also similar, then by the construction of $G'$, we can
obtain that $v_2$ and $v_3$ are similar in $G$. Since $G$ does not
contain $G_4$ and $G_5$ (see Figure \ref{fig1}) as an induced
subgraph, the set $N=N_G(v_1)\setminus\{u,v_2\}$ induces either a
complete graph or the disjoint union of two complete graphs.

{\bf Case 1.} $N$ induces a complete graph.

It is easy to observe that $G-v_2$ is connected and has no similar
vertices. Then $G-v_2$ has a $5$-leaf root $T'$. If there exists a
vertex $w\in N$ such that $d_{T'}(v_1, w)=5$, let $v_1abcdw$ be the
path in $T'$ joining vertices $v_1$ and $w$. We can construct a
$5$-leaf root $T$ of $G$ from $T'$ by attaching a path with length
$2$ ending at a new leaf $v_2$ to the vertex $c$.

In the following, we assume $d_{T'}(w,v_1)\leq 4$ for all $w\in N$.
If there exists a vertex $w\in N$ such that $d_{T'}(v_1, w)=4$, let
$v_1abcw$ be the path in $T'$ joining vertices $v_1$ and $w$. We can
construct a $5$-leaf root $T$ of $G$ from $T'$ by attaching a
pendent path with length $2$ ending at a new leaf $v_2$ to the
vertex $b$. Similarly, if $d_{T'}(w,v_1)=3$ for all vertices $w\in
N$, then attaching a path of length $2$ ending at the new leaf $v_2$
to the neighbor of the leaf $v_1$ in $T'$ yields a $5$-leaf root of
$G$. These contradictions complete the proof of {\bf Case 1}.

{\bf Case 2.} $N$ induces the disjoint union of two complete graphs
$N_1$ and $N_2$.

Since $G$ is a $4$-leaf power, $G$ is chordal by Proposition
\ref{pro1}. Thus the graph $G-\{u,v_1,v_2\}$ has exactly two
components with vertex sets $U_1$ and $U_2$ satisfying that
$N_1\subseteq U_1$ and $N_2\subseteq U_2$. Note that for $i=1,2$,
$|N_i|\geq 1$.

{\bf Subcase 2.1.} $|U_1|\geq 2$, $|U_2|\geq 2$.
\begin{figure}[ht]
\psfrag{u}{$u$}\psfrag{v1}{$v_1$}\psfrag{v2}{$v_2$}
\psfrag{a1}{$a_1$}\psfrag{a2}{$a_2$}\psfrag{b1}{$b_1$}\psfrag{b2}{$b_2$}
\psfrag{w1}{$w_1$}\psfrag{wt}{$w_t$}\psfrag{c1}{$c_1$}\psfrag{cs}{$c_s$}
\psfrag{F1}{$F_1$}\psfrag{F2}{$F_2$}
\begin{center}
\includegraphics[width=10cm]{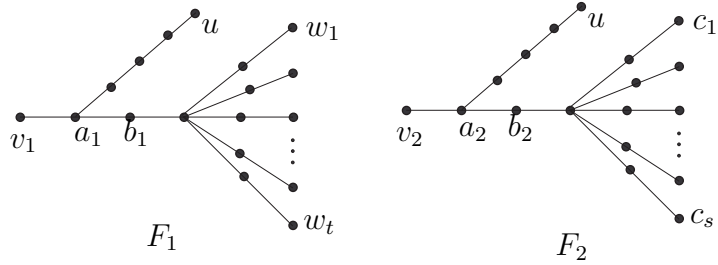}
\end{center}
\caption{$5$-leaf roots of graphs $H_1$ and $H_2$ in Subcase 2.1.}
\label{fig3}
\end{figure}

Since $|U_i|\geq 2$ for $i=1,2$ and $G$ has no pair of similar
vertices, we have $U_i\setminus N_i\neq \emptyset$. For any vertex
$x_i\in U_i\setminus N_i$ ($i=1,2$), $|N_G(x_i)\cap N_i|\leq 1$. If
$|N_i|=1$ for $i=1,2$, we are done. Suppose $|N_i|\geq 2$, let
$N_1=\{w_1,w_2,\ldots, w_t\}$ and $N_2=\{c_1,c_2,\ldots, c_s\}$.
Suppose for some $i$, there is a vertex $x_i\in U_i\setminus N_i$
and $|N_G(x_i)\cap N_i|\geq 2$. Then $G$ contains $G_8$ (see Figure
\ref{fig1}) as a subgraph, a contradiction. Let $H_1=G-U_2-v_2$ and
$H_2=G-U_1-v_1+uv_2$. For $i=1,2$, $H_i$ is a $5$-leaf power without
pairs of similar vertices by the choice of the graph $G$. It is
obvious that $H_1, H_2\in \mathcal{H}$. By Lemma \ref{lem1}, there
must be a $5$-leaf root $T_i$ of $H_i$ containing $T_0$ as its
subtree, denote by $F_i$ the subtree of $T_i$ which is isomorphic to
$T_0$, as shown in Figure \ref{fig3}. We can construct a new tree
$T$ from $T_1$ and $T_2$ by deleting the pendent path containing $u$
in $T_2$, adding an edge $b_1b_2$ and contracting edges $a_1b_1$ and
$a_2b_2$. It is easy to verify that $T$ is a $5$-leaf root of $G$, a
contradiction.

\begin{figure}[ht]
\psfrag{u}{$u$}\psfrag{v1}{$v_1$}
\psfrag{a2}{$a_2$}\psfrag{b2}{$b_2$}
\psfrag{c1}{$c_1$}\psfrag{cs}{$c_s$} \psfrag{F2}{$F_2$}
\begin{center}
\includegraphics[width=5cm]{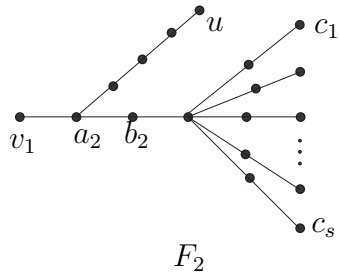}
\end{center}
\caption{$5$-leaf root of graph $H_2$ in Subcase 2.2.} \label{fig4}
\end{figure}

{\bf Subcase 2.2.} $|U_1|=1$, $|U_2|\geq 2$.

Suppose $U_1=\{w_1\}$, $N_2=\{c_1,c_2,\ldots, c_s\}$. By similar
discussions as Subcase 2.1, we know that $U_2\setminus N_2\neq
\emptyset$ and for any vertex $x\in U_2\setminus N_2$, $|N_G(x)\cap
N_2|\leq 1$. Let $H_2=G-w_1-v_2$. Since $H_2$ has no pairs of
similar vertices, $H_2\in \mathcal{H}$. By Lemma \ref{lem1}, $H_2$
has a $5$-leaf root $T_2$ containing $T_0$ as its subtree, denote by
$F_2$ the subtree of $T_2$ which is isomorphic to $T_0$, as shown in
Figure \ref{fig4}. Now we also can construct a $5$-leaf root $T$ of
$G$ by attaching a pendent path with length $2$ ending at a new leaf
$v_2$ to the vertex $b_2$ and attaching a pendent path with length
$2$ ending at a leaf $w_1$ to the vertex $a_2$, a contradiction.
 \qed

From Proposition \ref{pro2} (i) and Theorem \ref{thm1}, we conclude
the following corollary.
\begin{cor}
If $G$ is a $4$-leaf power with $L(G)\neq \emptyset$, then $G$ is
also a $k$-leaf power for $k\geq 4$.
\end{cor}

\end{document}